\documentclass[12pt, amssymb, amsmath]{article}

\begin{document}

\renewcommand{\thesection}{\arabic{section}}
\renewcommand{\thesubsection}{\thesection.\arabic{subsection}}


\newtheorem{lem}{Lemma}[section]
\newtheorem{propo}[lem]{Proposition}
\newtheorem{theo}[lem]{Theorem}
\newtheorem{rema}[lem]{Remark}
\newtheorem{remas}[lem]{Remarks}
\newtheorem{coro}[lem]{Corollary}
\newtheorem{defin}[lem]{Definition}
\newtheorem{hypo}[lem]{Hypoth\`ese}
\newtheorem{exem}[lem]{Exemple}
\newtheorem{conj}[lem]{Conjecture}

\newcommand{\real}{{\bf R}}
\newcommand{\reall}{{\bf R}^{p}}
\newcommand{\realll}{{\bf R}^{k}}
\newcommand{\mreal}{M\times{\bf R}^{p}}
\newcommand{\mreall}{M\times{\bf R}^{2k}}
\newcommand{\ri}{\rightarrow}
\newcommand{\tast}{T^{\ast}M}
\newcommand{\cd}{D_{\bullet}}
\newcommand{\cdk}{D_{\bullet-k}}
\newcommand{\cfx}{C_{\bullet}(f,\xi )}
\newcommand{\cald}{{\cal D}(M)}
\newcommand{\cm}{C_{\bullet}(M)}
\newcommand{\cmk}{C_{\bullet-k}(M)}
\newcommand{\cov}{\widetilde{M}}
\newcommand{\calh}{{\cal FH}}
\newcommand{\calf}{{\cal F}}
\newcommand{\calc}{{\cal C}}
\newcommand{\calx}{\cal X}
\newcommand{\caly}{{\cal Y}}
\newcommand{\gol}{\, \, \, \, \, \, \, \, }
\newcommand{\hatf}{\hat{f}}
\newcommand{\hatx}{\hat{x}}
\newcommand{\hatv}{\hat{v}}
\newcommand{\mn}{M^{n}}
\newcommand{\barm}{\bar{M}}
\newcommand{\bfs}{{\bf S}}
\newcommand{\piu}{\pi_{1}}
\newcommand{\pii}{\pi_{i}(M)}

\begin{center}
{\large\bf ON THE  HOMOTOPY  OF FINITE CW-COMPLEXES WITH POLYCYCLIC FUNDAMENTAL
GROUP}\\

\vspace{.3in}

\noindent Mihai DAMIAN\\
\noindent Universit\'e Louis Pasteur\\
IRMA, 7, rue Ren\'e Descartes,\\
67 084 STRASBOURG\\
e-mail : damian@math.u-strasbg.fr\\

 june 2006

\end{center}
\vspace{.4in}

\noindent{\bf Abstract\, :}\,  Let $X$ be a finite 
connected CW-complex of dimension $q$. If its fundamental group
$\piu(X)$ is polycyclic of Hirsch number $h>q$, we show that  at least one
homotopy group $\pi_{i}(X)$  is not finitely generated. If $h=q$ or $h=q-1$ the
same conclusion holds unless $X$ is an Eilenberg-McLane space $K(\piu(X),1)$. 
\vspace{.2in}

\section{Introduction}

  Let $X$ be a finite connected CW-complex of dimension $q$. 
 Consider the homotopy groups 
$\pi_{i}(X)=[\bfs^{i},X]$, for $i\geq 2$. If all these 
(Abelian) groups are finitely generated we say that the homotopy 
of $X$ is finitely generated.  For simply connected complexes, 
a celebrated theorem of Serre \cite{Ser} asserts that

\begin{theo}\label{serre} If $\piu(X)=1$ the homotopy of $X$ is finitely generated.
\end{theo}

 On the other hand  simple examples as 
 $X\, =\, \bfs^{1}\vee\bfs^{2}$ show that the homotopy of $X$ is not always finitely
generated. 

  One can then ask if there is a general negative statement asserting 
that the homotopy of $X$ is not finitely generated under some hypothesis 
on computable invariants of $X$, such as its  fundamental group. This is 
the aim of the present paper.  The statements we will prove have the 
following form :

\begin{quote}
If the fundamental group of $X$ satisfies to the conditions (C) (which  
depend on the dimension $q$), then the homotopy of $X$ is not finitely 
generated unless $X$ is an Eilenberg-Mc.Lane space $K(\pi_{1}(X))$. 
\end{quote}

  If $X$ is an Eilenberg-Mc Lane space, the cohomological dimension of
$\piu(X)$ is less or equal than $q\, =\, dim(X)$. So, we may add $cd(\pi_{1}(X))>q$  to 
the hypothesis $(C)$ in order to have the desired conclusion on the homotopy 
of $X$. For instance  when $\piu(X)$ has non trivial torsion elements then
$cd(\piu(X))$ is infinite and the above conclusion is valid.\\
 
\subsection{Statement of the results}

  Before stating the main theorem we remind the following :

\noindent {\bf Definition.} A group $G$ is called polycyclic if it admits a
series $$1=G_{0}<G_{1}<\cdots <G_{k}=G$$ 
 with cyclic factors $G_{i+1}/G_{i}$. When all the factor groups are infinite cyclic
we call $G$ a poly-${\bf Z}$ group.\\

K.A. Hirsch proved in \cite{Hir1} that the number of infinite 
cyclic factors in such a series is  an invariant of $G$. It 
is called Hirsch number of $G$ and it is denoted by $h(G)$. 

Here are our main statements :

\begin{theo}\label{main} Suppose that $X$ is a finite connected CW complex of
dimension $q$ and that $\piu(X)$ is polycyclic. Then \\
a) If $h(\piu(X))\, >\, q$, the homotopy of $X$ is not finitely generated. \\
b) If $h(\piu(X))\, \in\, \{q-1,q\}$ then the homotopy of $X$ is not finitely
generated unless $X$ is a $K(\piu(X),1)$. In particular, when $\piu(X)$ has torsion
then the conclusion of a) holds. 
\end{theo}

 Note that, by a result of J-P. Serre (see \ref{Ser} below),
 when the homotopy of $X$ is not finitely
generated at least one of the groups $\pi_{2}(X), \pi_{3}(X),\ldots, \pi_{q}(X)$ is
not finitely generated.  When $X$ is a manifold we can improve this result :

\begin{theo}\label{mainman} Let $\mn$ be a closed manifold. Let $r\, =\, 
\max\left\{\left[\frac{n}{2}\right],3\right\}$. Suppose that $\piu(M)$ is polycyclic.
Then, if $h(\piu(M))\geq n-1$, the group $\pii$ is not finitely generated for some
$i\leq r$ unless $h(\piu(M))=n$ and $M$ is a $K(\piu,1)$. 

Moreover, if $n\geq 6$ we can replace the assertion "M is a $K(\piu,1)$" by the
stronger one "the universal cover of $M$ is diffeomorphic to $\real^{n}$".
\end{theo}

When $n=4$, we are able to prove the above statement for $r=2$, namely~:

\begin{theo}\label{fourmain} Let $M^{4}$ be a closed connected manifold with
polycyclic fundamental group. If $h(\piu(M))\geq 3$ then $\pi_{2}(M)$ is not finitely
generated unless $h(\pi_{1}(M))=4$ and $M$ is a $K(\piu,1)$.\end{theo}

In the next theorem  we weaken the hypothesis on the fundamental group but
we suppose in addition that $\chi(X)\neq 0$. Recall first the finiteness 
properties of a group, which were introduced by C.T.C. Wall in \cite{Wa} :\\

\noindent{\bf Definition.} Let $r\geq 1$ be an integer. A groupe $G$ is of 
type $\calf_{r}$
 if there is an Eilenberg-McLane space $K(G,1)$ whose $r$-skeleton has a finite 
number of cells.  Equivalently a group is of type $\calf_{r}$ if it acts freely, 
properly, cellularly and cocompactly on an $(r-1)$-connected cell complex.

 A group is of type $\calf_{\infty}$ if it is of type $\calf_{r}$ for any integer
$r>0$. \\

\noindent{\bf Remark.} A group $G$ is of type $\calf_{1}$ if and only if 
it is finitely 
generated. $G$ is of type $\calf_{2}$ if and only if it is finitely presented.\\

\begin{theo}\label{chi} Let $X^{q}$ be a finite connected 
CW-complex of dimension $q$ with
fundamental group of type $\calf_{q+1}$. Suppose that there is a non vanishing morphism
$u:\piu(X)\ri{\bf Z}$ such that $Ker(u)$ is of type $\calf_{q+1}$. Then, if $\chi(X)
\neq 0$, the homotopy of $X$ is not finitely generated. 

If $X$ is a manifold with non-zero Euler characteristic we obtain the same
holds for $\piu(X)$ and $Ker(u)$ of type $\calf_{\left[\frac{q}{2}\right]}$.\end{theo}

The hypothesis on $u$ may be reformulated in terms of the Bieri-Renz invariants of
$\piu(X)$ as follows :
$$(\ast)\gol \Sigma^{q+1}(\piu)\, \cap\, -\Sigma^{q+1}(\piu)\, \neq\, \emptyset.$$
The Bieri-Renz invariants $\Sigma^{i}(G)$ of a group $G$ 
are open subsets of the unit sphere of
$\real^{rk(G)}$. We recall the definition and the properties 
of these invariants in Section 4. 

In the next subsection we show that polycyclic fundamental groups satisfy the
hypothesis of \ref{chi} 
 and
we give other examples of groups for which the condition $(\ast)$ is fulfilled.

\subsection{Comments on the results}

  An Abelian group is obviously polycyclic of Hirsch number equal to its rank.  
It is also of type $\calf_{\infty}$. For Abelian fundamental  groups 
the results of the
previous subsection were proved by the author 
in \cite{dami}. Here are other remarks about these statements. \\

  \noindent {\bf Remarks.} 1. The lower bound $q-1$ 
for $h(\piu)$ in the hypothesis of \ref{main} is optimal. Indeed, the complex $X= 
{\bf T}^{n-2}\times{\bfs}^{2}$  fulfills 
the conditions of the hypothesis of \ref{main} except for  $h(\piu(X))=q-2$.
 It is a consequence 
of the theorem \ref{serre} of J.P.Serre, quoted above, that $X$ does 
not satisfy to the conclusion of \ref{main}. \\
\\ 2. If $G$ is a polycyclic group then $[G,G]$ is obviously polycyclic and
therefore for every morphism $u:G\ri{\bf Z}$, $Ker(u)$ has the same
property. This implies that the hypothesis on the fundamental group 
 of \ref{chi} is fulfilled when
$\piu(X)$ is polycyclic. Indeed we have :

\begin{propo}\label{polycalf} Polycyclic groups are of type $\calf_{\infty}$.
\end{propo}

\noindent\underline{Proof}

The proposition is an immediate corollary of the following :

\begin{lem}\label{polyf} Consider an exact sequence of groups $$1\ri K\ri G\ri Q\ri
1.$$
If $K$ and $H$ are of type $\calf_{\infty}$ then so is $G$. 
\end{lem} 

  To prove the above lemma one may use the following results :\\
a) $G$ is $\calf_{\infty}$ iff $H_{\ast}(G,\cdot)$ commutes with direct products.
This was proved by R. Bieri and B. Eckmann in \cite{BE}.\\
b) The Hochshild-Serre spectral secquence \cite{hs} which satisfies
$$E_{pq}^{2}=H_{p}(Q, H_{q}(K,R))$$ and converges towards $H_{\ast}(G,R)$. (The
complete statement \ref{hochschild} is given in Section 4). 

It is clear how a) and b) imply \ref{polyf} proving thus \ref{polycalf}. \hfill
$\diamond$\\
\\
3. Other examples of groups satisfying the hypothesis of \ref{chi} were constructed
by M. Bestvina, N. Brady in \cite{BB} and then in a more general statement by J.
Meier, H. Meinert and L. Van Wyk in \cite{MMW}. Let us describe them briefly : \\
\\
{\bf Definition.} A flag complex $L$ is a finite 
simplicial complex with the property that any
collection of $q+1$ mutually adjacent vertices span a $q$-simplex in $L$. The {\it
right angled Artin group} associated to $L$ is the group $G_{L}$ spanned by the
vertices $v_{1}, \ldots, v_{s}$ of $L$ with relations $[v_{i},v_{j}]=1$ whenever
$v_{i}$ and $v_{j}$ are adjacent in $L$. \\

 Note that the group $G_{L}$ admits a finite Eilenberg-McLane space \cite{BB} 
and therefore it is
of type $\calf_{\infty}$.  We have

\begin{theo}\label{BB} (Bestvina, Brady, Meier, Meinert, van Wyk) Let $L$ and $G_{L}$
as above and consider a morphism $u:G_{L}\ri {\bf Z}$ such that $u(v_{i})\neq 0$ for
each generator $v_{i}$. Then, if $L$ is $(r-1)$ connected, $Ker(u)$ is of type
$\calf_{r}$.\end{theo}

Considering appropriate flag complexes $L$ (remark that the barycentric subdivision of
any complex is a flag complex) we may find a lot of examples of groups
which fulfil the hypothesis of \ref{chi}. 

 \subsection{Idea of the proof}

  The general idea of the proof of \ref{mainman} is the following : we show 
first that a manifold $M$ as in \ref{mainman} whose dimension is greater or 
equal to $6$ and whose homotopy groups $\pii$ are of finite type for $i\leq r$ 
admits a  fibration over the circle.   
This is the main difficulty of the proof. To overcome it, we will use   
Novikov homology theory and some of its applications on  Morse functions 
$f:M\ri\bfs^{1}$. 

Remark that if the manifold  $F^{n-1}$ 
is a  fiber, the inclusion $j:F\hookrightarrow M$ induces  
isomorphisms in $\pi_{i}$ for $i\geq 2$. The lift of $j$ induces a homotopy 
equivalence between the universal covers of $F$ and $M$.\\

  Now take a manifold $M$, as in \ref{mainman} and consider the product 
$M\times\bfs^{3}$ in order to fulfill the condition on the dimension. 
Supposing that the homotopy groups $\pii$ are of finite type for $i\leq r$, 
apply the previous argument to this product and get a fiber $F$ as above. 
Then  check that $F$ still satisfies the hypothesis of \ref{mainman}.
 If its dimension is still greater or equal than $6$ we apply the same 
argument to $F$.  By succesive iterations we find a manifold $F_{0}$ of low 
dimension whose universal cover has the same homotopy type as the one of 
$M\times\bfs^{3}$. In particular the homology groups of 
$\widetilde{F}_{0}$ and $\widetilde{M}\times\bfs^{3}$ are isomorphic. 
By comparing them, we infer that the universal cover of $M$ is acyclic, 
which means that $M$ is Eilenberg-McLane. So, either the homotopy of 
$M$ is not finitely generated, or $M$ is Eilenberg-Mc.Lane, as in the 
statement of \ref{mainman}.

To prove \ref{main} we embed $X$ into an Euclidian space and thicken it to a
manifold $W$, with boundary $M$. Then we apply the above argument to $M$.

In the hypothesis of \ref{chi}, supposing that the homotopy of $X$ is finitely
generated, we are only able to prove that {\it the Novikov homology} $H_{\ast}(M,u)$
of the corresponding manifold $M$ vanishes. But this implies that $\chi(M)=\chi(X)=0$,
yielding a contradiction.\\

The paper is organized as follows. In Section 2 we state  
the result \ref{fibre} on the fibration over the circle which 
was roughly sketched above. Supposing \ref{fibre} true, we show how it 
can be succesively applied in order to prove \ref{mainman}, \ref{fourmain} 
 and \ref{main}. In Section 3 
we recall the definition and some useful properties of the Novikov homology. 
We also recall some basic facts about Morse theory of circle-valued functions 
and point out the relation between the (vanishing of the) Novikov homology and 
the existence of a  fibration over the circle. In Section 4 we prove 
\ref{fibre} and \ref{chi}. We will use the Bieri-Renz criterion which we recall 
in the Subsection 4.3, devoted to the Bieri-Renz invariants.

\section{Iterated fibrations}

  Our main result \ref{main} is a consequence of the following :

\begin{theo}\label{fibre}  Let $\mn$ a closed manifold of dimension $n\geq 6$. 
Suppose that $\piu$ is of type $\calf_{\left[\frac{n}{2}\right]}$. Suppose also 
that $\pii$ are of finite type for $i\leq \left[\frac{n}{2}\right]$. 

Suppose that there is a non zero cohomology class  
$u\in H^{1}(M;\real)\, \approx\, Hom(\piu(M),\real)$  such that $Ker(u)$ 
is of type $\calf_{\left[\frac{n}{2}\right]}$. Suppose also that the Whitehead
group $Wh(\piu(M))$ vanishes. 

Then there is a fibration $f:M\ri {\bf S}^{1}$ such that $[f^{\ast}d\theta]=u$. 

When  $M$ of  arbitrary dimension $n$ and $Ker(u)$ is of type $\calf_{r}$, 
where $r\, =\, \max\left\{\left[\frac{n}{2}\right],2\right\}$ the same 
conclusion holds if we replace $M$ by  $M\times {\bf
S}^{p}$ for all $p\geq 6-n$.
\end{theo}
\noindent{\bf Remarks.} 1. The Whitehead group is defined as follows :
$$Wh(\pi)\, =\, \frac{GL{\bf Z}[\pi]}{\left[GL{\bf Z}[\pi],GL{\bf Z}[\pi]\right],\{\pm
g|g\in\pi\}},$$
where $GL{\bf Z}[\pi]=\lim_{\ri}GL_{m}{\bf Z}[\pi]$. \\
2. T. Farell and W. Hsiang proved in \cite{FH} that $Wh(\pi)$ vanishes when $\pi$ is
poly-${\bf Z}$. The Whitehead group also vanishes when $\pi=G_{L}$ is one of the
examples of Bestvina and Brady. This result was proved by B. Hu \cite{Hu}. 

It is  
conjectured that $Wh(\pi)=0$ for any torsion-free group $\pi$.\\
3. For $\piu=\bf Z$, Theorem \ref{fibre} was proved by W. Browder and J. Levine in
\cite{BL}. \\

 The proof of \ref{fibre} will be given in  Section 4. Let us now show  how this 
theorem implies
our main results \ref{mainman}, \ref{fourmain} and \ref{main}. \\

\noindent \underline{Proof of \ref{fibre} $\Longrightarrow$ \ref{mainman}}\\

 Without restricting the generality of our statements, we may suppose that 
$n\geq 3$.  We begin with the statement of the following result, due 
to K. Hirsch (\cite{Hir2}, Theorem 2).

\begin{theo}\label{hirsch} Let $G$ be a polycyclic group. There exists a normal
subgroup $N$ of finite index in $G$ which is poly-${\bf Z}$ and such that $h(N)=h(G)$. 
\end{theo}

  Now, let $\piu^{0}\leq\piu$ a subgroup as in \ref{hirsch} and consider the
associated finite cover $M_{0}\ri M$. Let   $$1=G_{0}<G_{1}<\cdots <G_{k}=\piu^{0}$$
be a series with infinite cyclic factor groups : we have therefore $h(\piu)=k$, so by
hypothesis $k\geq n-1$.  Denote by $u_{1}$ the projection
$G_{k}\ri G_{k}/G_{k-1}\approx {\bf Z}$.\\
 
   Suppose  that $\pii$ is finitely generated for 
$i\leq r = \max\left\{\left[\frac{n}{2}\right],3\right\}$. 

Consider first the case $n\geq 6$.  Using 
\ref{polycalf} we find that $u_{1}$ fulfills the hypothesis of 
\ref{fibre}. We apply
 \ref{fibre} and we get a fibration $M_{0}\ri{\bf S}^{1}$. Let $F_{1}$ be a fibre of
this fibration.   So $F_{1}$  satisfies :\\
1.  $\pi_{1}(F_{1})=Ker(u_{1})=G_{k-1}$. \\
2. $\pi_{i}(F_{1})\approx\pi_{i}(M)$ $\forall\,  i \geq 2$.\\

 We wish to apply \ref{fibre} to $F_{1}$. We use the cohomology class 
$u_{2}:G_{k-1}\ri G_{k-1}/G_{k-2}$. 
 Its kernel is $G_{k-2}$, so, using again 
\ref{polycalf}, $Ker(u_{2})$ is of type $\calf_{\left[\frac{n}{2}\right]}$ (actually
$\calf_{\infty}$, since polycyclic).
 The hypothesis on the higher homotopy groups of \ref{fibre} is fulfilled by 
$F_{1}$ (because of the condition 2. above), therefore we may apply 
\ref{fibre} to $F_{1}$ if its dimension is no less than $6$. 

  We get thus a closed connected manifold $F_{2}\subset F_{1}$ whose higher 
homotopy groups are those of $M$ and whose fundamental group is 
$G_{k-2}$. If its dimension is greater
or equal than $6$ we may again apply \ref{fibre} to the couple
$(F_{2},u_{3}:G_{k-2}\ri G_{k-2}/G_{k-3})$. 

  By iterating this argument we get a secquence
$$ (1)\gol  F_{n-5}\hookrightarrow
F_{n-6}\hookrightarrow\cdots\hookrightarrow F_{1}\hookrightarrow M_{0}$$
 
such that for $j=1,\ldots, n-1$ :\\
1. $\pi_{1}(F_{j})= G_{k-j}$\\
2. $dim(F_{j}) =n-j$.\\
3. The inclusion $F_{j}\hookrightarrow M$ induces isomorphisms in $\pi_{i}$ for
$i\geq 2$. \\

 The manifold $F_{n-5}$ is of dimension $5$, so it  does not verify the 
dimension hypothesis of \ref{fibre}. Its fundamental group is $G_{k-n+5}$.
In order to continue to apply 
\ref{fibre}, we consider the
product $F\, = \, F_{n-5}\times {\bf S}^{3}$. As above, we have the 
cohomology class $u_{k-n+4}:G_{k-n+5}\ri G_{k-n+5}/G_{k-n+4}$
 which is non zero and has a (poly-${\bf Z}$) kernel of 
type $\calf_{\left[\frac{n}{2}\right]}$.  By \ref{fibre}, we get a 
fibration of $F$ over ${\bf S}^{1}$. Its fiber $K_{0}$  is a closed
connected manifold of dimension $7$ with finitely generated $\pi_{i}$ for
$i=1,\ldots,\, r= \max\left\{\left[\frac{n}{2}\right],3\right\}$ 
(since $\bfs^{3}$ has
the same property by \ref{serre}). Its fundamental group is 
$G_{k-n+4}$.
 Since $[7/2]=3$, we may apply \ref{fibre} to
$K_{0}$ and then again to its submanifold  $K_{1}$ given by 
\ref{fibre},  to obtain a sequence  as above, where the
maps are inclusions of fibers of fibrations over the circle
 which therefore induce isomorphisms at the level of 
$\pi_{i}$ for $i\geq 2$ :
$$ K_{2}\hookrightarrow K_{1}\hookrightarrow K_{0} \hookrightarrow F.$$

It follows that the universal covers of $K_{2}$ and $F$ are homotopically
equivalent, in particular
$$(2)\gol
H_{\ast}(\widetilde{K}_{2})\, \approx\, H_{\ast}(\widetilde{F}_{n-5}\times\bfs^{3}).$$

Now $K_{2}$ is a closed connected $5$-dimensional manifold whose fundamental 
group is
$G=G_{k-n+2}$. Since $k\geq n-1$, $G$ is an infinite group.  It follows that 
$H_{5}(\widetilde{K}_{2})\, =\, 0$
(see \cite{Bredon}, p.346) ; We have therefore that $H_{i}(\widetilde{K}_{2})$ 
vanishes for $i>4$.

 Using
 the Kunneth formula we infer from (2) that $H_{i}(\widetilde{F}_{n-5})$ vanishes for
$i>0$. Therefore $\widetilde{F}_{n-5}$ is contractible and, using the sequence (1)
$\widetilde{M_{0}}=\widetilde{M}$ is contractible, too. So $M$
 is an Eilenberg-Mc Lane space $K(\pi_{1}(M),1)$. The cohomological dimension of
$\piu(M)$ is therefore equal to $n$.  Now we use the following well-known result (for a
proof, see \cite{Grun},
Lemma 8, p. 154) : 

\begin{propo}\label{grunberg}

If $G$ is a polycyclic group then $cd(G)$ is finite if and only if $G$ is torsion
free. In this case $cd(G)=h(G)$. 
\end{propo}

 Applying \ref{grunberg} we get that $M$ cannot be a $K(\piu,1)$ unless
$h(\piu(M))=n$. If this relation is not valid, we get a contradiction and
therefore  the homotopy of $M$ cannot be
finitely generated, completing the proof of \ref{mainman} for
$n\geq 6$. If it is and if $M$ is a $K(\piu,1)$, let us show that its universal cover
is $\real^{n}$. Since according to \ref{fibre} $M$ fibers over ${\bf S}^{1}$, we have 
$$\widetilde{M}\, \approx\, \widetilde{F}\times \real.$$
 To finish the proof we just have to apply a celebrated theorem of J. Stallings
\cite{Sta} which asserts :

\begin{theo}
If a manifold $P^{n\geq 5}$ is a product of two open non-trivial contractible
manifolds than $P$ is diffeomorphic to $\real^{n}$. 
\end{theo}

  \hfill $\diamond$\\

\noindent\underline{Proof of \ref{fibre} $\Longrightarrow$ \ref{fourmain} (case n=4)}\\

  Suppose that $\pi_{2}(M)$ is finitely generated. As in the proof of \ref{mainman},
use \ref{hirsch} and consider a finite cover $M_{0}$ with poly-${\bf Z}$ fundamental
group $\piu^{0}$ of Hirsch number greater than $3$. We apply first 
\ref{fibre} to $M_{0}\times\bfs^{2}$ and to the morphism $u_{1}:G_{k}\ri
G_{k}/G_{k-1}$.
We get a manifold $F$ of dimension $5$ which is a fiber of a fibration $M_{0}\ri {\bf
S}^{1}$. Then we apply again \ref{fibre} to the product 
$F\times \bfs^{2}$ and the cohomology class $u_{2}: G_{k-1}\ri G_{k-1}/G_{k-2}$ 
and we get a $6$-dimensional submanifold $K \hookrightarrow F\times\bfs^{2}$.  
We have a secquence 
$$ K\hookrightarrow F\times\bfs^{2}\hookrightarrow 
M \times\bfs^{2}\times\bfs^{2},$$
which induces homotopic equivalences at the level of the universal covers. 
In particular 
$$H_{\ast}(\widetilde{K})\, \sim\, H_{\ast}(\widetilde{M}\times\bfs^{2}
\times\bfs^{2}).$$
 Now, as above, since $k\geq 3$, the fundamental group $G_{k-2}$ of 
 $K$ is infinite, so the homology of 
$\widetilde{K}$ vanishes in degrees $i>5$. This implies that the homology of 
$\widetilde{M}$ is zero, so  $M$ is a $K(\piu,1)$. We conclude using \ref{grunberg}.
\hfill $\diamond$\\

\noindent\underline{Proof of \ref{fibre} $\Longrightarrow$ \ref{mainman} 
 in the cases $n = 3, \, 5$}\\

 Suppose that $\pi_{2}(M)$ and $\pi_{3}(M)$ are finitely generated. 
We apply $(n~-~2~)$ times Theorem \ref{fibre} to the $8$-dimensional manifold 
$M\times\bfs^{8-n}$. We get  a secquence :
$$K_{n-2}\hookrightarrow K_{n-3}\cdots\hookrightarrow M\times\bfs^{8-n},$$
as above. The manifold $K_{n-2}$ is of dimension $8-n+2$ and has an infinite
 fundamental group. So, $H_{i}(\widetilde{K}_{n-2})=0$ for $i\geq 8-n+2$. 
Since this homology is isomorphic to $H_{\ast}(\widetilde{M}\times\bfs^{8-n})$, 
it follows that $\widetilde{M}$ is acyclic, and therefore that $M$ is 
an Eilenberg-Mc.Lane. As above, we use \ref{grunberg} to finish the proof. \hfill
$\diamond$\\

\noindent\underline{Proof of \ref{fibre} $\Longrightarrow$ \ref{main}}\\

 Assume that the homotopy of $X$ is finitely generated. 

 a)  Suppose  that $h(\piu(X))= k$ for some integer $k>q$. Using \ref{hirsch},
consider   a finite cover $X_{0}$  of $X$ such that $\piu(X_{0})$ is poly-$\bf Z$ and
$h(\piu(X_{0}))=k$. 

 Embed $X_{0}$ in an Euclidian space $\real^{2q+r+1}$ for some $r\geq 0$ which will
be fixed later in the proof. 
Let W be a tubular neighbourhood of $X_{0}$ and denote by
$M^{2q+r}$
the smooth manifold $\partial W$. Since $M$ is a deformation retract of
$W\setminus X_{0}$, using a
general position argument we get  isomorphisms between  $\pi_{i}(W)$ and
 $\pi_{i}(M)$ for $i\leq q+r-1$. In particular, for $r>1$ the inclusion
$M\hookrightarrow W$ induces an isomorphism at the level of fundamental groups and,
since $X_{0}$ is a retract of $W$, the higher homotopy groups $\pi_{i}(M)$ 
are finitely generated   for $i \leq q+r-1$.
 If $r\geq 1$  then  $q+r-1\geq 
\left[\frac{2q+r}{2}\right]$ and the hypothesis on the higher homotopy groups 
of \ref{fibre} is
satisfied by the manifold $M$. Let $$1=G_{0}<G_{1}<\cdots <G_{k}=\piu(X_{0})$$
be a series with infinite cyclic factor groups.   For $r$ large enough
 ($r\geq k-2q+5$), we may apply $k$ times \ref{fibre}, as 
in the proof of \ref{mainman}, 
 and get a sequence of closed
connected manifolds :

$$ (3)\gol  F_{k}\hookrightarrow
F_{k-1}\hookrightarrow\cdots\hookrightarrow F_{1}\hookrightarrow M$$
 
such that for $j=0,\ldots, k$ :\\
1. $\pi_{1}(F_{j})=G_{k-j}$. \\
2. $dim(F_{j}) =2q+r-j$.\\
3. The inclusion $F_{j}\hookrightarrow M$ induces isomorphisms in $\pi_{i}$ for
$i\geq 2$. \\

The constant $r\geq 1$ must be chosen large enough to have :\\
a) $dim(F_{j})\geq 5\, \, 
 \forall j\, $ : needed for the dimension hypothesis in \ref{fibre} (this
means $ 2q+r-k \geq 5$). \\
b) $r+q-1\geq 
\left[\frac{2q+r}{2}\right]$, (which is true for $r\geq 1$)
 to insure the hypothesis on the higher homotopy groups in
\ref{fibre}, as we explained above.\\
c) $r>k-q$. \\

It follows that the first manifold  $F_{k}$ in the sequence (1) is closed and simply
connected, of dimension $2q+r-k>q$. Using the property 3 we find that the  
 composition
 $$(4)\gol  F_{k}\hookrightarrow M\hookrightarrow W \ri X_{0}$$ induces
isomorphisms at the level of the higher homotopy groups $\pi_{i}$ for $i=2, \ldots,
q+r-1$. By Whitehead's theorem, the induced application between the universal covers
of $F_{k}$ and $X_{0}$ is an isomorphism at the level of the $i^{th}$ homology group
for $i\leq q+r-1$. In particular, since $q+r-1\geq 2q+r-k$ (since by hypothesis 
 $k>q$), we have 
$$(5) \gol H_{2q+r-k}(\widetilde{F}_{k})\, \approx\, H_{2q+r-k}(\widetilde{X}_{0}).$$
 But $F_{k}$ is simply connected of dimension $2q+r-k$, so the left side of (5) is 
 $\bf Z$. On the other hand $2q+r-k\, >\, q=dim(X_{0})$, so the right side of (5)
vanishes. This contradicts the initial assumption on the finite generation of the
homotopy of $X$   and the 
statement a) of the theorem
 is proved. \\

b) Suppose now that $k=h(\piu(X))\, \in\, \{q-1,q\}$. As above we get a closed
simply connected manifold $F_{k}$ of dimension $2q+r-k$ (i.e. $q+r+1$ or $q+r$) such 
that the application $H_{i}(F_{k})\, \ri \, H_{i}(\widetilde{X_{0}})$ is an
isomorphism for $i\leq q+r-1$. In particular, the homology 
$H_{i}(F_{k};{\bf Z})$ vanishes
for $i= q+1, \ldots, q+r-1$. 

By the univesal coefficients theorem we obtain that the
cohomology $H^{i}(F_{k};{\bf Z})$ is zero  for $i=q+2, \ldots, q+r-1$, so
Poincar\'e duality implies that $H_{i}(F_{k};{\bf Z})$ also vanishes for 
for $i=q-k+1,\ldots, q+r-k-2$. Now $q-k+1\in\{1,2\}$ and $F_{k}$ is simply connected :
We infer after choosing  $r$  sufficiently large 
that the integer homology of $F_{k}$ vanishes
in the degrees $i\leq q$ which means that $\widetilde{X}_{0}$ is contractible.
Therefore $X_{0}$ and $X$ are Eilenberg-Mac Lane spaces. In particular the 
cohomological dimension of $\pi_{1}(X)$ cannot exceed $q$,
 which implies  that $\pi_{1}(X)$ is torsion free. 

This completes the proof. \hfill $\diamond$

\section{ Novikov homology and fibrations over the circle}

 In the  preceeding section we showed that our main results \ref{main},  
\ref{mainman} and \ref{fourmain}
are consequences of  \ref{fibre}.  Let $M$ be a closed 
connected manifold and $u\in H^{1}(M;{\bf Z})\approx Hom(\piu,{\bf Z})$. 
The aim of the present 
section is to introduce the Novikov homology $H_{\ast}(M,u)$ and to describe the
situations when the vanishing of $H_{\ast}(M,u)$ implies the existence of a fibration
$f:M\ri{\bf S}^{1}$ i.e. the conclusion of \ref{fibre}. Then in the Section 4 we
 prove that the hypothesis of \ref{fibre} implies the vanishing of the Novikov
homology $H_{\ast}(M,u)$.  

\subsection{Novikov homology}

Let $u\in H^{1}(M;\real)$.
Denote by $\Lambda$ the ring   $ {\bf Z}[\pi _{1}(M)]$ and by $ \hat{\Lambda}$ the
ring of formal series ${\bf Z}[[\pi_{1}(M)]]$.
Consider a $\calc^{1}$-triangulation of $M$
which we lift  to the universal
cover $\widetilde{M}$. We get a $\Lambda$-free complex $C_{\bullet}(M)$  spanned by
(fixed lifts of) the cells of the triangulation of $M$.  \\

We define now {\it the completed ring} $\Lambda_{u}$  :
$$\Lambda_{u}:=\, \left\{\lambda\, =\, \sum n_{i} g_{i}\, \in\,
\hat{\Lambda}\, \, |\, \, g_{i}\in \pi_{1}(M), \, \, n_{i}\in {\bf Z}, \, \, \,
u(g_{i})\, \rightarrow\,  +\infty\right\}$$
The convergence to $+\, \infty$ means
 here that for all $A\, >\, 0$,  $u(g_{i})\, <\, A$ only for a finite
number of $g_{i}$ which appear with a non-zero coefficient in  the sum  $\lambda$.\\

\noindent{\bf Remark} Let $\lambda\, =\, 1+\sum n_{i} g_{i}$ 
where $u(g_{i})>0$ for all $i$.
Then $\lambda$ is invertible in $\Lambda_{u}$. Indeed, if we 
denote by $ \lambda_{0}\, =\,
\sum n_{i} g_{i}$ then it is easy to check that $\sum_{k\geq 0}(-\lambda_{0})^{k}$ 
is an element of
$\Lambda_{u}$ and it is obvious that it is the inverse of $\lambda$.\\

\noindent{\bf Definition} Let $C_{\bullet}(M,u)$ be the  
$\Lambda_{u}$-free complex
$\Lambda_{u}\otimes_{\Lambda}C_{\bullet}(M)$. The Novikov homology $H_{\ast}(M,u)$ is
the homology of the complex
$C_{\bullet}(M,u)$.\\

A purely algebraic  consequence of the previous definition is the following version of
the universal coefficients theorem (\cite{Gode}, p.102, Th 5.5.1) :

\begin{theo}\label{gode} There is a spectral sequence $E_{pq}^{r}$ which converges to
$H_{\ast}(M,u)$ and such that
$$E_{pq}^{2}\, =\, Tor_{p}^{\Lambda}( H_{q}(\widetilde{M}),\Lambda_{u}).$$
\end{theo}

We will use this result in Section 4 to prove that in the hypothesis of 
\ref{mainman}, the Novikov
homology associated to some class vanishes. 

\subsection{Morse-Novikov theory}

We recall in this subsection the relation between Novikov homology and 
closed one forms. In dimension $n\geq 6$, when the Novikov homology vanishes, 
some hypothesis on $\piu$ stated below imply the existence of a nowhere  
vansihing closed one form on $M$. It is well-known (see \cite{Ti}) that 
the existence of such a form is equivalent to the existence of a 
fibration of $M$ over $\bfs^{1}$. 

Let $\alpha$ be a closed generic one form in the class $u$.  Let $\xi$ be the
gradient of $\alpha$ with respect to some generic metric on $M$. For every critical
point $c$ of $\alpha$ we fix a point $\tilde{c}$ above $c$ in the universal cover
$\widetilde{M}$.  We can define then a
complex $C_{\bullet}(\alpha,\xi)$ spanned by the zeros of $\alpha$ :
the incidence number $[d,c]$ for two zeros
of consecutive indices is the (possibly infinite) sum
 $\sum n_{i} g_{i}$ where $n_{i}$ is the algebraic number of
flow lines which join $c$ and $d$ and which are covered by a path in $\widetilde{M}$
 joining
$g_{i}\tilde{c}$ and $\tilde{d}$. It turns out that this incidence number belongs to
$\Lambda_{u}$, so $C_{\bullet}(\alpha,\xi)$ is actually a $\Lambda_{u}$-free complex.

 The fundamental property of the Novikov homology is that it is isomorphic to the
homology of the complex $C_{\bullet}(\alpha,\xi)$
above for any couple $(\alpha,\xi)$.

By comparing the complexes $C_{\bullet}(\alpha,\xi)$ and $C_{\bullet}(-\alpha,-\xi)$
 we get the following duality property
(see Prop.2.8 in \cite{Da} and 2.30 in \cite{Lat}) :

\begin{theo}\label{dual} Let $\mn$ be a closed connected manifold, $u\, \in\,
H^{1}(M;\real)$ and let $l$ be an integer. If $H_{i}(M,-u)=0$ for $i\leq l$, then
$H_{i}(M,u)=0$ for $i\geq n-l$.
\end{theo}

If the form $\alpha$ has no zeroes then $C_{\bullet}(\alpha,\xi)$ vanishes and
therefore we have  $H_{\ast}(M,[\alpha])~=~0$. Conversely,
one can ask if the vanishing of
$H_{\ast}(M,u)$ implies the existence of a nowhere vanishing closed $1$-form belonging
to the class $u\in H^{1}(M)$. For $n\geq 6$
this problem was independently solved by F. Latour
 \cite{Lat} and A. Pajitnov \cite{Paj}, \cite{Paj'}.
The statement is (\cite{Lat}, Th.1') :
 \begin{theo}\label{Latour} For $dim(M)\geq 6$ the following set of
conditions is equivalent to the existence of a nowhere vanishing closed
 $1$-form in  $u\in H^{1}(M,{\bf Z})$ :

1. Vanishing Novikov homology  $H_{\ast}(M,u)$.

2. Vanishing Whitehead torsion $\tau(M,u)\, \in \, Wh(M,u)$.

3. Finitely presented $Ker(u)\subset\pi_{1}(M)$.
\end{theo}

  \noindent{\bf Remarks} \\
1. The definition of the generalized Whitehead group $Wh(M,u)$
 and of the Whitehead torsion is given in \cite{Lat}. \\
2. In the statement of \cite{Paj} the first two conditions are replaced by
:

1'. $C_{\bullet}(M,u)$ is simply equivalent to zero.

Actually, one can show (see \cite{Mau}) that 1' is equivalent to "1 and 2". \\
3. In earlier works on the subject as those of F.T. Farell \cite{Fa1} and 
L. Siebenmann
\cite{Sie} the algebraic conditions which are equivalent to the existence of a 
nowhere
vanishing closed $1$-form in a rational cohomology class $u$ were stated in the
hypothesis that the infinite cyclic cover $P$ associated to $u$ is finitely dominated.
($P\ri M$ is defined to be the pull-back of the universal covering
$\real\ri {\bf S}^{1}$ defined by a function $f:M\ri{\bf S}^{1}$ such that
$[f^{\ast}d\theta]\, =\, u$). The  relation between the finite
domination of $P$ and vanishing of the Novikov homology  was first established by
 A. Ranicki in \cite{Ra2}, \cite{Ra3}. A Ranicki (\cite{Ra1}, Chap.14)
 established a relation  between the
Whitehead "fibering obstruction" from \cite{Fa1} and \cite{Sie} and the condition 2
above. Then A. Pajitnov and A. Ranicki proved in \cite{Para}  the following  
\begin{theo}\label{para}
 If $G$ is a group with $Wh(G)=0$ and $u:G\ri{\bf Z}$ is a morphism, then $Wh(G,u)=0$. 
\end{theo}

  It follows that under the hypothesis of \ref{fibre} the condition 2 of \ref{Latour}
is always satisfied. 

In order to prove \ref{chi} we use

\begin{propo}\label{novichi} If $H_{\ast}(M,u)=0$ for some class $u$ then $\chi(M)=0$.
\end{propo}

\noindent\underline{Proof}

 The complex $C_{\bullet}(M,u)$ is acyclic.  Like in
\cite{Mau} one can then show that this  complex
 is  simply equivalent to a complex of the form
$$0\ri\Lambda_{u}^{p}\stackrel{\partial}{\ri}\Lambda_{u}^{p}\ri 0.$$
This means that the second complex is isomorphic to the first after adding or
cancelling afinite number of trivial summands  
$$0\ri\Lambda_{u}\stackrel{Id}{\ri}\Lambda_{u}\ri 0.$$
In particular $$\chi(M)=\chi(C_{\bullet}(M,u))=0.$$
\hfill $\diamond$

\section{Novikov homology and finiteness properties of groups}

In this section we will achieve the proof of our results \ref{main} ... \ref{chi}.
 Up to now we have 
proved that \ref{main}, \ref{mainman} and \ref{fourmain} are implied by 
\ref{fibre} (Section 2) and that the
conclusion of \ref{fibre} is valid if the three conditions of \ref{Latour} are
fulfilled. It remains to show 
that the hypothesis of \ref{fibre} implies these three conditions. It is clearly 
the case for the third one, since in the hypothesis of \ref{fibre}, $Ker(u)$ 
is of type 
$\calf_{\left[\frac{n}{2}\right]}$, in particular of type $\calf_{2}$, i.e. 
finitely presented. The second condition is also satisfied, according to \ref{para}.

 The main point of our proof is the vanishing of the Novikov 
homology $H_{\ast}(M,u)$ under the assumption of the hypothesis of \ref{fibre}.

  Recall that, by \ref{gode}, there is a spectral secquence $E_{pq}^{r}$  which 
converges to $H_{\ast}(M,u)$ and whose term $E_{pq}^{2}$ is equal to 
$Tor_{p}^{\Lambda}( H_{q}(\widetilde{M}),\Lambda_{u})$. Recall also that, 
by the duality result \ref{dual} we have the implication :
$$ H_{i}(M,\pm u)=0 \, \, \forall i \leq\left[\frac{n}{2}\right]\, \, 
\, \, \Rightarrow\, \, \, H_{i}(M,u)=0\, \forall i.$$

 It suffices therefore to prove the following statement :

 \begin{propo}\label{important} Let $\mn$ a closed manifold. Suppose that 
$\piu(M)$ is of type $\calf_{\left[\frac{n}{2}\right]}$. Suppose also that 
$\pii$ are of finite type for $i\leq \left[\frac{n}{2}\right]$. 

Suppose that there is a non zero cohomology class  $u\in H^{1}(M;\real)\, 
\approx\, Hom(\piu(M),\real)$  such that $Ker(u)$ is of 
type $\calf_{\left[\frac{n}{2}\right]}$. Then, for all integers 
$0\, \leq\, p,q\, \leq\, \left[\frac{n}{2}\right]$ we have
$$Tor_{p}^{\Lambda}( H_{q}(\widetilde{M}),\Lambda_{\pm u})\, =\, 0.$$
\end{propo}  

Note that this result is purely algebraic. In order to proof it we use 
some facts about

\subsection{Hurewicz-type morphisms}

Recall that the classical Hurewicz's theorem asserts that for $q\geq 2$ 
the cannonical
morphism $I_{q}:\pi_{q}(M)\ri H_{q}(M)$ is an isomorphism provided that $M$ is
$(q-1)$-connected.

In \cite{Ser} J-P. Serre generalized this theorem (see also \cite{Sp}, p. 504) : For
some "admissible" classes of groups $\calc$, he showed that,
if   $X$ is simply connected such that  $\pi_{i}(X)\, \in\, \calc$ for
$i=1,\ldots,q-1$, where $q\geq 2$,
 then $I_{q}$ is an isomorphisme mod $\calc$ : This means
that  $Ker(I_{q})$
and $Coker(I_{q})$ are in $\calc$.  

The class of finitely generated Abelian groups is such an admissible class. In
particular we have :

\begin{theo}\label{Ser} Let $X$ be a simply connected space. Then $\pi_{i}(X)$ is
finitely generated for $i\leq q$ iff $ H_{i}(X)$ is finitely generated for $i\leq q$.

In particular any closed, simply connected CW-complex has finitely generated homotopy
groups (which is \ref{serre}).
\end{theo}

  By applying this theorem to $\widetilde{M}$ we may replace the hypothesis 
on  $\pii$ by the analogue hypothesis on $H_{i}(\widetilde{M})$. From now on we 
will suppose that for all $i\leq \left[\frac{n}{2}\right]$, $H_{i}(\widetilde{M})
={\bf Z}^{r_{i}}\oplus T_{i}$, where $T_{i}$ is a torsion finitely generated 
$\bf Z$-module.   The  proof of \ref{important} relies on the following statement :

\begin{propo}\label{decisiv} Let $\pi$ be a group of type $\calf_{p}$ for some 
positive integer $p$ and $u:\pi\ri\bf Z$ a morphism whose kernel is of type 
$\calf_{p}$. Suppose that ${\bf Z}^{r}$ is a $\pi$-module. Denote by $\Lambda$ 
the ring ${\bf Z}[\pi]$ and by $\Lambda_{u}$ the completed ring, as above. Let 
$\pi_{0}\leq \pi$ be a normal subgroup of finite index and $\Lambda_{0}$ the 
corresponding group ring ${\bf Z}[\pi_{0}]$. 
Then for $i\leq p$ we have $$ Tor_{i}^{\Lambda_{0}}( {\bf Z}^{r},
\Lambda_{ u})\, =\, 0.$$
\end{propo}

\noindent{\bf Remark}  For $\pi_{0}=\pi$ and $r=1$ (and therefore 
for arbitrary $r$ and trivial action of $\pi$ on ${\bf Z}^{r}$)
 the statement above was proved by J-C. Sikorav in this thesis 
\cite{Sikt}. \\

We postpone the proof of \ref{decisiv} to Subsection 4.4. We now show :\\

\subsection{Proof of \ref{decisiv}\, $\Longrightarrow$\,  
\ref{important}}

 Fix $p,q\leq \left[\frac{n}{2}\right]$. For $g\in\piu(M)$, denote by 
$\phi_{g}$ the automorphism of $H_{q}(\widetilde{M})
={\bf Z}^{r_{q}}\oplus T_{q}$, given by the action of $\piu(M)$. We have :
$$\phi_{g}\, =\, \left(\begin{array}{ccc}
a_{g}&\, &0\\
b_{g}&\, & c_{g}
\end{array}\right),$$
where $a_{g}:{\bf Z}^{r_{q}}\ri{\bf Z}^{r_{q}}$, $b_{g}:{\bf Z}^{r_{q}}\ri T_{q}$ 
and $c_{g}:T_{q}\ri T_{q}$. Note that $a_{g}$ and $c_{g}$ are automorphisms 
(of inverses $a_{g^{-1}}$, resp. $c_{g^{-1}}$). Let 
$$\pi_{0}\, =\, \{\, g\in\piu\, |\, c_{g}\, =\, Id\}.$$
  Since $T_{q}$ is finite, there is only a finite number of automorphisms 
$c:T_{q}\ri T_{q}$. Therefore $\pi_{0}$ is a normal subgroup of $\piu(M)$ of 
finite index (It is the kernel of the morphism $\piu(M)\, \ri\,  Aut(T_{q})$ 
defined by $g\mapsto c_{g}$). We use the following :

\begin{lem}\label{index}  Let $G$ a group of type $\calf_{p}$ and 
$G_{0}\leq G$ a normal subgroup of finite index. Then $G_{0}$ is of 
type $\calf_{p}$. 
\end{lem}

  The proof of \ref{index} is obvious : If $Q$ is a $K(G,1)$ with finite $p$-skeleton, 
then the finite cover of $Q$ corresponding to $G_{0}\leq G$ will be a $K(G_{0},1)$ 
with finite $p$-skeleton. \\

 Now let $u:\piu(M)\ri \real$, as in the statement of \ref{important} and let 
$u_{0}\, =\, u|_{\pi_{0}}$. Obviously, $Ker(u_{0})$ has finite index in 
$Ker(u)$ so, using \ref{index}, both $\pi_{0}$ and $Ker(u_{0})$  are of 
type $\calf_{p}$.  

  Consider now the short exact sequence :
$$0\ri T_{q}\stackrel{(0,Id)}{\longrightarrow} {\bf Z}^{r_{q}}\oplus 
T_{q}\stackrel{Id\oplus 0}{\longrightarrow} 
 {\bf Z}^{r_{q}}\ri 0.$$

One immediately checks that this is an exact sequence of $\pi_{0}$-modules 
(where the action of $\pi_{0}$ on ${\bf Z}^{r_{q}}$ is $x\mapsto a_{g}(x)$).  
Now consider $\Lambda_{0}={\bf Z}[\pi_{0}]$, and view the completed ring 
$\Lambda_{u}$ as a $\Lambda_{0}$-module. The tensor product of $\Lambda_{u}$ 
and the exact sequence above yields a   long exact sequence :
$$(1)\gol  \cdots\ri Tor_{p}^{\Lambda_{0}}( T_{q},\Lambda_{ u})\ri 
Tor_{p}^{\Lambda_{0}}( {\bf Z}^{r_{q}}\oplus T_{q},\Lambda_{ u})\ri 
Tor_{p}^{\Lambda_{0}}( {\bf Z}^{r_{q}},\Lambda_{ u})\ri\cdots$$

  As  a consequence of \ref{decisiv}, the right term in the sequence 
above vanishes. In order to prove that the left term is zero, 
we consider an exact sequence of the form 
$$0\ri {\bf Z}^{s}\ri{\bf Z}^{m}\ri T_{q}\ri 0,$$
which is viewed as a sequence of $\pi_{0}$-modules with trivial actions 
(recall that, by construction, $\pi_{0}$ acts trivially on $T_{q}$). 
The corresponding long exact sequence given by the tensor product with 
$\Lambda_{u}$ writes :
$$ (2)\gol \cdots\ri Tor_{p}^{\Lambda_{0}}( {\bf Z}^{m},\Lambda_{ u})\ri
 Tor_{p}^{\Lambda_{0}}( T_{q},\Lambda_{ u})\ri Tor_{p-1}^{\Lambda_{0}}
( {\bf Z}^{s},\Lambda_{ u})\ri\cdots.$$

  Applying again \ref{decisiv} we infer 
$Tor_{p}^{\Lambda_{0}}( T_{q},\Lambda_{ u})=0$, therefore the 
middle term in $(1)$ vanishes. \\

  We have thus  established that for each $p,q\leq \left[\frac{n}{2}\right]$
$$(3) \gol Tor_{p}^{\Lambda_{0}}( H_{q}(\widetilde{M}),\Lambda_{ u})\, =\, 0,$$
which is the assertion required in \ref{important} ... with $\Lambda_{0}$ 
instead of $\Lambda$. To complete the proof we need the following results : 

\begin{propo}\label{brown} Let $G$ be a groupe and 
let $L$ be a ${\bf Z}[G]$ right module and $N$ be a ${\bf Z}[G]$ left module. 
Assume that $N$ is ${\bf Z}$-torsion-free. Then 
$$ Tor^{{\bf Z}[G]}_{\ast}(L,N)\, \approx\,  H_{\ast}(G, L\otimes_{\bf Z}N), 
$$ where $G$ acts diagonally on $L\otimes_{\bf Z}N$ : 
$g(x\otimes y)=xg^{-1}\otimes gy$. 
\end{propo}

This proposition is proved in \cite{brown} (prop. 2.2, p. 61).

\begin{theo}\label{hochschild} For any group extension 
$$1\ri K\ri G \ri Q\ri 1$$ and any $G$-module $R$ there 
is a spectral sequence of the form $$E_{ij}^{2}\, =\, H_{i}(Q,H_{j}(K,R))$$ 
which converges to $H_{i+j}(G,R)$.
\end{theo}

  This theorem is due to G. Hochschild and J-P. Serre \cite{hs} 
(see also \cite{brown}, p. 171).  

We apply first \ref{brown} and infer from $(3)$ that for each 
$p,q\leq \left[\frac{n}{2}\right]$ we have :
$$ (4)\gol H_{p}(\pi_{0}, H_{q}(\widetilde{M})\otimes_{\bf Z}\Lambda_{u})
\, =\, 0.$$
Note that the hypothesis of \ref{brown} is fulfilled by $N=\Lambda_{u}$.

 We fix $q$ and denote by $R$ the $\piu(M)$-module 
$H_{q}(\widetilde{M})\otimes_{\bf Z}\Lambda_{u}$ (for the diagonal action). 
Then we apply \ref{hochschild} to the extension 
$$1\ri\pi_{0}\ri \piu(M)\ri\piu(M)/\pi_{0}\ri 1$$ and to the module $R$ ;
We find using $(4)$  that  $E_{ij}^{2}=0$ for all 
$j\leq \left[\frac{n}{2}\right]$ and for all $i\in\bf N$. 
According to \ref{hochschild}, this implies that 
$$ (5)\gol H_{i}(\piu(M),R)=0\, \, \, \forall i\leq \left[\frac{n}{2}\right].$$

  Finally, we apply once again \ref{brown} and we get 
that for all $i,q\leq \left[\frac{n}{2}\right]$ we have~:
$$Tor_{i}^{\Lambda}( H_{q}(\widetilde{M}),\Lambda_{u})\, =\, 0.$$
 
The analogous relation for $\Lambda_{-u}$ instead of $\Lambda_{u}$ 
can be established in the same way, so \ref{important} follows.\hfill $\diamond$\\

  So we only have to proof \ref{decisiv} to complete the proof of \ref{main},
\ref{mainman} and \ref{fourmain}. 
The proof involves some facts about Bieri-Renz invariants. We recall the 
definition and some properties of these invariants in the subsection below :

\subsection{Bieri-Renz invariants. Proof of \ref{decisiv}} 

Let $G$ be a group of type $\calf_{m}$. We call two non zero homomorphisms 
$u,v:G\ri \real$ equivalent if $u=\lambda v$ for some positive $\lambda\in \real$. 
We denote by $S(G)$ the quotient $Hom(G,\real)/\approx$ 
(which is a $(rk(G)-1)$-dimensional sphere).  The Bieri-Renz 
invariants $\Sigma^{i}(G)$ and $\Sigma^{i}(G,{\bf Z})$, 
defined for $i=1,\ldots,m$  are open subsets of $S(G)$. 
They were introduced by R. Bieri, W Neuman and R. Strebel in \cite{BNS} for $i=1$ and
by  R. Bieri and B. Renz in \cite{BR} for $i\geq 2$.  

  These invariants are defined as follows. Let $X$ be a $K(G,1)$ 
which is a complex with finite $m$-squeleton and let $u:G\ri \real $
 be a non zero homomorphism. Then there exists an {\it equivariant height function} 
$f:\widetilde{X}\ri \real$, i.e. a function which satisfies $f(gx)=f(x)+u(g)$. 
(If $X$ is a manifold then $f$ is a primitive of the pullback of some $1$-form 
in the class $u$). It can be shown that the difference of two such functions 
is bounded. Consider the maximal subcomplex $X_{f}$ of $\widetilde{X}$ whose 
image by $f$ is contained in $[0, +\infty[$. \\
\\
{\bf Definition.}  \, {\it Let $u\in S(G)$ and $i\in\{1,\ldots,m\}$. 
Then $u$ belongs to $\Sigma^{i}(G)$ (resp. to $\Sigma^{i}(G,{\bf Z})$) 
if for some couple $(X,f)$ as above the subcomplex $X_{f}$ is $(i-1)$-connected 
(resp. (i-1)-acyclic).} \\

 Bieri and Renz show that the definition does not depend on $(X,f)$.  
The following properties of the Bieri-Renz invariants are obvious from 
the definition :\\
a)  $\Sigma^{i}(G)\subset\Sigma^{i}(G,{\bf Z})$.\\
b)   $\Sigma^{1}(G)=\Sigma^{1}(G,{\bf Z})$.\\
c) For $i\geq 2$, $\Sigma^{i}(G)=\Sigma^{i}(G,{\bf Z})\cap\Sigma^{2}(G)$.\\

The most striking application of these invariants is stated in the following~:

\begin{theo}\label{renz} Let $N\subset G$ be a normal subgroup with 
Abelian quotient $G/N$. denote by $S(G,N)$ the subset of $S(G)$ defined 
by $\{u\in S(G)\, |\, u|_{N}=0\}$. Then, for any $i=1,\ldots,m$, we have the
equivalence :\\
 $N$ is of type $\calf_{i}$ 
iff $S(G,N)\subset \Sigma^{i}(G)$. \end{theo}

 Note that $\{ u,-u\}\subset S(G,Ker(u))$ ; When $Im(u)$ is cyclic it 
is easy to show that  these two sets coincide. Note also that 
$S(G,[G,G])=S(G)$. So, one immediately infers the following corollary :

\begin{coro}\label{cororenz}
i) Let $u\in S(G)$ and $i$ an integer as above. If 
$Ker(u)$ is $\calf_{i}$ then $\pm u\in\Sigma^{i}(G)$. 
The converse is valid if $Im(u)$ is cyclic. \\
ii) $\Sigma^{i}(G)=S(G)$ if and only if $[G,G]$ is of type $\calf_{i}$.\end{coro}

  In the sequel we give an algebraic description of the 
invariants $\Sigma^{i}(G,{\bf Z})$ following \cite{BR}, Section 4.
 Fix a non zero homomorphism 
$u:G\ri \real$. Let $F$ be a finitely generated free ${\bf Z}[G]$-module, 
and $\{e_{i}\}_{i=1,\ldots, k}$ a  basis of $F$. We define an application
 $v:F\ri\real$ as follows : $v$ is defined arbitraily on  the elements $e_{i}$ ;
 Then for any $g\in G$ we put $v(ge_{i})=v(e_{i})+u(g)$. Finally, if
 $\lambda\in F$ writes $\lambda=\sum_{i,j}n_{ij}g_{j}e_{i}$ in the basis 
$\{e_{i}\}$, we define $$v(\lambda)\, =\, inf\{v(g_{j}e_{i})\, |\, n_{ij}\neq 0\},$$
and $v(0)=+\infty$. Following Bieri and Renz we call $v$ a 
{\it valuation extending $u$}. Now suppose that $G$ is of type $\calf_{m}$ and let 
$$P_{m}\stackrel{\partial_{m}}{\ri} P_{m-1}\stackrel{\partial_{m-1}}{\ri}
\cdots \stackrel{\partial_{1}}{\ri} P_{0}\stackrel{\partial_{0}}{\ri}{\bf Z}\ri 0$$
be a ${\bf Z}[G]$-free, finitely generated resolution (which exists since 
there exists a $K(G,1)$ with finite $m$-skeleton). 
Define as above $v_{i}:P_{i}\ri\real$ which are valuations extending $u$. 
We may suppose in addition that for any $i=1,\ldots, m$ and for any 
$x\in P_{i}$ we have $$(1) \gol v_{i}(x)\leq v_{i-1}(\partial_{i}(x)).$$ 
Indeed, one easily sees that it suffices to check the relation above on the 
basis $\{e_{i}\}$ of $P_{i}$. We can construct $v$ inductively, by choosing 
$v(e_{i})$ sufficiently negative in order to satisfy the inegality $(1)$. 

  Bieri and Renz proved the following theorem (\cite{BR}, theorem 4.1) : \\

\begin{theo}\label{renz2} Let $P_{\bullet}\ri{\bf Z}$ be a ${\bf Z}[G]$-free,
 finitely generated resolution of length $m$ and let $v:P_{\bullet}\ri\real$ be a 
valuation extending $u$ satisying (1). Then $u\in\Sigma^{m}(G,{\bf Z})$ if
 and only if there exists a chain  endomorphism $\Phi:P_{\bullet}\ri P_{\bullet}$ 
which lifts the identity of ${\bf Z}$ and which satisfies the property 
$$(2)\gol v(\Phi(x))>v(x)\, \, \, \, \forall\, x\in P_{\bullet}.$$
\end{theo}

We will reformulate this theorem as follows. It is obvious that one 
 can construct the valuation 
$v:P_{\bullet}\ri\real$ such that it satisfies an additional feature~: for all
 $j=1,\ldots,m$ $v$  is constant  on the set $\{e^{j}_{i}\}$ of the basis 
elements of $P_{j}$. Denote this constant by $\nu_{j}$. 
Then, the valuation $v_{j}:P_{j}\ri \real$ is given by $$v_{j}(\lambda)\, =\, \nu_{j}
+inf\{u(g_{k})\, |\, n_{ik}^{j}\neq 0\},$$
where  $\lambda =\sum_{i,k} n_{ik}^{j}g_{k}e^{j}_{i}$, 
$g_{k}\in G$,  $n_{ik}^{j}\in \bf
Z$.
 Thus, for an endomorphism $\Phi$ given by \ref{renz2}, 
if $$\Phi(e^{j}_{i})\, =\, \sum_{i,k} n^{j}_{ik}g^{j}_{ik} e^{j}_{k}\, ,$$
(where $g_{ik}^{j}\in G$ and $n_{ik}^{j}\in \bf Z$), the inegality $(2)$ 
implies that $u(g_{ik}^{j})>0$ for all elements of $G$ appearing with non 
zero coefficient $n_{ik}^{j}$. \\

 We call an element $\lambda=\sum_{i}n_{i}g_{i}$ of $ {\bf Z}[G]$ 
$u$-positive if $u(g_{i})>0$ for any $g_{i}$ which has a non zero 
coefficient $n_{i}$ in the writing of $\lambda$. We call a matrix 
$A\in{\cal M}_{k}({\bf Z}[G])$ $u$-positive if all its entries are 
$u$-positive. Taking into account the preceeding remarks, 
theorem \ref{renz2} can be stated as follows :

\begin{theo}\label{renz3}  Let $P_{\bullet}\ri{\bf Z}$ a resolution of 
length $m$ which is ${\bf Z}[G]$-free and finitely generated.  
Fix a basis for each $P_{j}$ for all $j=1,\ldots, m$.  Then 
$u\in\Sigma^{m}(G,{\bf Z})$ if and only if there exists a chain  
endomorphism $\Phi:P_{\bullet}\ri P_{\bullet}$ which lifts the 
identity of ${\bf Z}$ and such that for all $j=1,\ldots m$ the 
matrices of $\Phi_{j}:P_{i}\ri P_{j}$ in the fixed basis are $u$-positive. 
\end{theo}

Now we are able to complete the \\

\noindent\underline{Proof of \ref{decisiv}} 

 Applying \ref{brown} we infer that 
$ Tor_{i}^{\Lambda_{0}}( {\bf Z}^{r},\Lambda_{u})$ is isomorphic to 
$H_{i}(\pi_{0}, {\bf Z}^{r}\otimes_{\bf Z}\Lambda_{u})$, 
where the action of $\pi_{0}$ on ${\bf Z}^{r}\otimes_{\bf Z}\Lambda_{u}$ is
 given by $x\otimes\lambda\mapsto xg^{-1}\otimes g\lambda$. It suffices 
therefore to prove that the latter vanishes for $i\leq p$. 

  Let $$P_{p}\ri P_{p-1}\ri\cdots\ri P_{1}\ri P_{0}\ri{\bf Z}$$ be a 
${\bf Z}[\pi_{0}]$-free resolution which we may suppose finitely 
generated since $\pi_{0}$ is $\calf_{p}$. By definition we have 
$$(3)\gol H_{i}(\pi_{0}, {\bf Z}^{r}\otimes_{\bf Z}\Lambda_{u}))
\approx H_{i}(P_{\bullet}\otimes_{{\bf Z}[\pi_{0}]}({\bf Z}^{r}
\otimes_{\bf Z}\Lambda_{u})).$$

  We will prove that the right term of $(3)$ vanishes for all $i\leq p$. 
Fix a basis $\{e^{j}_{i}\}_{i}$ for each module $P_{j}$. Since 
$Ker(u|_{\pi_{0}})$ is of type $\calf_{p}$, it follows by \ref{cororenz}.i that 
$u\in\Sigma^{p}(\pi_{0},{\bf Z})$, so, applying \ref{renz3} we obtain an 
endomorphism $\Phi:P_{\bullet}\ri P_{\bullet}$ such that for all $j=1,\ldots, m$ 
 the matrix of $\Phi_{j}$ associated to the basis $\{e_{i}^{j}\}$ is $u$-positive.  

   Let $\Psi=\Phi\otimes Id: P_{\bullet}
\otimes_{{\bf Z}[\pi_{0}]}({\bf Z}^{r}\otimes_{\bf Z}
\Lambda_{u})\ri P_{\bullet}\otimes_{{\bf Z}[\pi_{0}]}
({\bf Z}^{r}\otimes_{\bf Z}\Lambda_{u})$. 

  The proof of \ref{decisiv} will be  complete if we prove the following :

\begin{lem}\label{finish} The homomorphism $Id-\Psi$ is invertible and 
it induces zero in homology.
\end{lem}

\noindent\underline{Proof}

Let us prove first that $Id-\Psi$ vanishes in homology. Since
 $\Phi:P_{\bullet}\ri P_{\bullet}$ lifts the identity of $\bf Z$,
 and $P_{\bullet}$ is free it is well known and easy to prove that 
there exists a homotopy $s:P_{\bullet}\ri P_{\bullet+1}$ between 
$\Phi$ and $Id$. It follows that $s\otimes Id$ is a homotopy between 
$\Psi$ and $Id$, so $Id-\Psi$ induces the zero morphism in homology.

   Now let us prove the first assertion of the lemma. Let 
$\{f_{1},f_{2},\ldots f_{r}\}$ be the canonical basis of ${\bf Z}^{r}$.
 We can see ${\bf Z}^{r}\otimes_{\bf Z}\Lambda_{u}$ as a right
 $\Lambda_{u}$-module endowed with the canonical structure : 
$(x\otimes\lambda)\mu=x\otimes\lambda\mu$. It is a free module of rank $r$
 and $\{f_{1}\otimes 1, f_{2}\otimes 1,\ldots,f_{r}\otimes 1\}$ is a basis 
for this module. 

Recall that we denoted  by $\Lambda_{0}$ the ring ${\bf Z}[\pi_{0}]$. 
For any $j=1,\ldots m$, the product
$P_{j}\otimes_{\Lambda_{0}}({\bf Z}^{r}\otimes_{\bf Z}\Lambda_{u})$ inherits of the
structure of right $\Lambda_{u}$-module described above.
 On the other hand, since $P_{j}$ is free
we have 
$$(4)\gol P_{j}\otimes_{\Lambda_{0}}({\bf Z}^{r}\otimes_{\bf Z}
\Lambda_{u})\, \approx  \, 
({\bf Z}^{r}\otimes_{\bf Z}\Lambda_{u})^{rk(P_{j})}\, \approx\, (\Lambda_{u})^{r\cdot
rk(P_{j})}.$$
If $\{e_{i}^{j}\}_{i=1,\ldots rk(P_{j})}$ is the given basis of $P_{j}$ then 
$$(5)\gol \{e_{i}^{j}\otimes_{\Lambda_{0}}(f_{s}\otimes 1)\}_{i=1,\ldots, rk(P_{j}), \,
s=1,\ldots r}$$ will be a basis for $P_{j}\otimes_{\Lambda_{0}}
({\bf Z}^{r}\otimes_{\bf Z}
\Lambda_{u})$. 

 It is easy to check that the isomorphisms $(4)$ preserve the right
$\Lambda_{u}$-module structure. Moreover, the differential
$\partial\otimes_{\Lambda_{0}} Id$ of the complex
$P_{\bullet}\otimes_{\Lambda_{0}}({\bf Z}^{r}\otimes_{\bf Z}\Lambda_{u})$ respects
this structure, therefore it is a complex free right $\Lambda_{u}$-modules.  \\

We claim that the matrices of $\Psi$ in the basis $(5)$ are $u$-positive. Indeed, for
fixed $j$, let 
$(\lambda_{ik})_{i,k=1,\ldots rk(P_{j})}$ be the matrix of $\Phi_{j}$ in the basis
$\{e^{j}_{i}\}$ ; This matrix is known as $u$-positive. We dropped down the
index $j$ from the coefficients of the matrix to simplify the notations. 
Denote by $\bar{\lambda}$
the image of an element $\lambda \in \Lambda_{0}$ under the endomorphism
$\Lambda_{0}\ri\Lambda_{0}$ induced by the involution of $\pi_{0}$ : $g\mapsto
g^{-1}$. The right action of
$\bar{\lambda}_{ik}$ on ${\bf Z}^{r}$, evaluated on the basis $\{f_{s}\}$ writes :
$$f_{s}\bar{\lambda}_{ik}\, =\, \sum_{l=1}^{r}n_{iks}^{l}f_{l}\, ,$$
for some integers $n_{iks}^{l}$. 

  We infer that $$\Psi_{j}(e_{i}^{j}\otimes_{\Lambda_{0}}(f_{s}\otimes
1))\, =\, \sum_{k}e_{k}^{j}\lambda_{ik}\otimes_{\Lambda_{0}}(f_{s}\otimes 1)\, =$$
$$=\,
\sum_{k}e_{k}^{j}\otimes_{\Lambda_{0}}(f_{s}\bar{\lambda}_{ik}\otimes 
\lambda_{ik})\, =\,
\sum_{k,l}e_{k}^{j}\otimes_{\Lambda_{0}}(n_{iks}^{l}f_{l}
\otimes 
\lambda_{ik})\, =\, $$
$$=\, \sum_{k,l}e_{k}^{j}\otimes_{\Lambda_{0}}(f_{l}
\otimes 
n_{iks}^{l}\lambda_{ik})\, =\, \sum_{k,l}[e_{k}^{j}\otimes_{\Lambda_{0}}(f_{l}
\otimes 1)]\, 
n_{iks}^{l}\lambda_{ik}\, ,$$
so the matrix of $\Psi_{j}$ is $u$-positive. 

Denote this matrix by $A_{j}$ and define another matrix $B_{j}$ by : $B_{j}=Id
+\sum_{k=1}^{+\infty}A_{j}^{k}$. As $A_{j}$ is $u$-positive it is easy to check that
the matrix $B_{j}$ belongs to ${\cal M}_{r\cdot rk(P_{j})}(\Lambda_{u})$. It therefore
defines using the basis $(5)$ an endomorphism of 
$P_{j}\otimes_{\Lambda_{0}}({\bf Z}^{r}\otimes_{\bf Z}
\Lambda_{u})$. It is actually an automorphism since obviously $(Id-A_{j})B_{j}\, =\,
Id$.   Finally, as $\Psi$ is a morphism of complexes,  the morphism $Id
+\Psi+\Psi^{2}+\cdots$ induced by $B_{j}$ will also commute with the differential. We
have finally got an automorphism of $P_{\bullet}\otimes_{\Lambda_{0}}
({\bf Z}^{r}\otimes_{\bf Z}
\Lambda_{u})$ whose inverse is $Id-\Psi$. This completes the proof of \ref{finish} and
hence the proof of \ref{decisiv} which implies our main theorems \ref{main},
\ref{mainman} and \ref{fourmain}. \hfill $\diamond$

 \subsection{Proof of \ref{chi}}

Assume that the homotopy of $X$ is finitely generated. 
By embedding $X$ in $\real^{2q+3}$ construct a manifold $M=\partial W$ of dimension
$2q+2$ as in the proof
of \ref{main}. We have $\chi(X)=\chi(W)=2\chi(M)$. 
By general position $\pi_{i}(M)\approx\pi_{i}(X)$ for $i\leq q+1$. These
groups are therefore finitely generated. The hypothesis of \ref{important} is
fulfilled and, by applying this result, we get $H_{\ast}(M,u)=0$. But according to
\ref{novichi} this implies $\chi(M)=0$, contradicting thus the hypothesis on
$\chi(X)$. 

If $X$ is a manifold, the result follows directly from \ref{important}, \ref{gode} and
\ref{novichi}.

The proof is finished. \hfill $\diamond$
\\
\\
{\bf Acknowledgements}  I am grateful to Thomas Delzant for our useful discussions on
the subject.

\vspace{.3in}

\noindent Mihai DAMIAN\\
\noindent Universit\'e Louis Pasteur\\
IRMA, 7, rue Ren\'e Descartes,\\
67 084 STRASBOURG\\
\\
e-mail : damian@math.u-strasbg.fr


\begin{thebibliography}{99}
 \bibitem{BB} M. Bestvina, N. Brady, \,  Morse theory and finiteness properties
of groups, \,  Invent. Math.  129 (1997), 445-470. 
\bibitem{BE} R. Bieri, B. Eckmann, \,  Finiteness properties of duality groups,
\, Comment. Math. Helv.  49 (1974), 74-83. 
\bibitem{BNS} R. Bieri, W. Neuman, R. Strebel, \,   A geometric invariant for
discrete groups, \, Invent. Math.  129 (1987), 451-477. 
\bibitem{BR} R. Bieri, B. Renz, \,  Valuations on free resolutions and higher
geometric invariants of groups, \, Comment. Math. Helv.  63 (1998), 464-497. 
\bibitem{Bredon} G. Bredon, \,  Topology and Geometry, \, GTM 139,  Springer
Verlag, New York, 1997.
\bibitem{BL} W. Browder, J. Levine \,   Fibering manifolds over the circle, \, 
Comment. Math. Helv.  40 (1966), 153-160. 
\bibitem{brown} K.S. Brown, \,  Cohomology of groups, \, GTM 87, Springer Verlag,
New York, 1982. 
\bibitem{Da} M. Damian, \,  Formes ferm\'ees non singuli\`eres et propri\'et\'es
de finitude des groupes\, Ann.Sci.Ec.Norm.Sup.  33 (2000) no.3, 301-320. 
\bibitem{dami} M. Damian, \,  On the higher homotopy groups of a finite
CW-complex, \, Topology and its Appl.  149 (2005), 273-284.
\bibitem{Fa1} F.T. Farell, \,  The obstruction to fibering a manifold over the
circle, \, Bull. Amer. Math. Soc.  73 (1967), 737-740. 
\bibitem{FH} F.T. Farell, W.C. Hsiang, \,  The Whitehead group of poly-(finite or
cyclic) groups, \, J. London Math. Soc. (2)  24 (1981), 308-324.
\bibitem{Grun}  K.W Gruenberg, \, 
 Cohomological Topics in Group Theory, L.N.M. 143,
Springer Verlag, New York, 1970. 
\bibitem{Gode} R. Godement, \,  Topologie alg\'ebrique et th\'eorie des
faisceaux, \, Hermann, Paris, 1958. 
\bibitem{Hir1} K.A. Hirsch, \,  On infinite soluble groups I, \, Proc. London Math.
Soc.  44 (1938), 53-60.
\bibitem{Hir2} K.A. Hirsch, \,  On infinite soluble groups IV, \, J. London Math.
Soc.  27 (1952), 81-85.
\bibitem{hs} G. Hochschild, J.P. Serre, \,  Cohomology of group extensions, \,
Trans. Amer. Math. Soc.  74 (1953), 110-134. 
\bibitem{Hu}  B. Hu \,  Whitehead groups of finite polyhedra with non-positive
curvature, \, Jour. of Diff. Geom.  38 3/1993, 501-517.
\bibitem{Lat} F. Latour, \,  Existence de 1-formes ferm\'ees non-singuli\`eres
dans une classe de cohomologie de de Rham, \, Publ. Math. IHES  80 (1994). 
\bibitem{Mau} S. Maumary, \,  Type simple d'homotopie, \, L.N.M.  48, Springer
Verlag, New York,  1967.
\bibitem{MMW} J. Meier, H. Meinert, L. Van Wyk, \,  Finiteness properties and
abelian quotients of graph groups, \, J. Math. Res. Lett.  3 (1996), 779-785. 
\bibitem{Mil} J. Milnor, \,  Whitehead torsion, \, Bull. Amer. 
Math. Soc.  72 (1966), \\358-426. 
\bibitem{Nov} S.P. Novikov, \,  Multivalued functions and functionals. An
analogue of the Morse theory, \, Soviet. Math. Dokl. vol  24, nr.  2/1981,
222-226.
\bibitem{Paj} A. Pajitnov, \,  Surgery on the Novikov complex, 
\, K-theory  10 (1996),
323-412.
\bibitem{Paj'} A. Pajitnov \,  Surgery on the Novikov complex, \, Rapport de
Recherche, Nantes (1993), \, http://193.52.98.6/~pajitnov.
\bibitem{Para} A. Pajitnov,  A. Ranicki, \,  The Whitehead group of the Novikov
ring, \, K-theory  21 (2000)  no.4, 325-365. 
\bibitem{Ra1} A. Ranicki, \,  High-dimensional knot theory, \, Springer (1998). 
\bibitem{Ra2} A. Ranicki, \,   Lower K- and L- theory, \, L.M.S. Lecture Notes,
Cambridge, 1992.
\bibitem{Ra3} A. Ranicki, \, 
 Finite domination and Novikov rings, \, Topology 
34 (1995), 619-632.
\bibitem{Sta} J. Stallings, \,  The piecewise linear structure of the Euclidean
space, \, Proc. Cambridge Phil. Soc.  58 (1962), 481-488.
\bibitem{Ser} J.P. Serre, \,  Groupes d'homotopie et
 classes de groupes ab\'eliens, \, 
Ann. of Math. (2)  58 (1953), 258--294. 
\bibitem{Sie} L. Siebenmann, \,  A total Whitehead obstruction to fibering over
the circle, \, Comment. Math. Helv.  45 (1970), 1-48. 
\bibitem{Sikt} J-C. Sikorav, \,   Homologie de Novikov associ\'ee \`a une classe
de cohomologie r\'eelle de d\'egr\'e un, \,  Th\`ese Orsay 1987.
\bibitem{Sp} E. Spanier, \,  Algebraic Topology, \, Mc Graw-Hill, New York,  1966.
\bibitem{Ti} D. Tischler, \,  On fibering certain foliated  manifolds over
$S^{1}$, \, Topology  9 (1970), 153-154.
\bibitem{Wa} C.T.C. Wall, \,  Finiteness conditions for CW-complexes, \, Ann.
of Math.  81 (1965), 56-69. 
 
\end{thebibliography}
\end{document}